\RequirePackage{ifpdf}
\ifpdf 
\documentclass[pdftex]{sigma}
\else
\documentclass{sigma}
\fi

\usepackage{bbm,eucal,oldgerm}

\newcommand{\ed}{{\rm d}}
\newcommand{\w}{{\mathchoice{\,{\scriptstyle\wedge}\,}{{\scriptstyle\wedge}}
      {{\scriptscriptstyle\wedge}}{{\scriptscriptstyle\wedge}}}}
\newcommand{\lhk}{\mathbin{\hbox{\vrule height1.4pt width4pt depth-1pt
             \vrule height4pt width0.4pt depth-1pt}}}

\renewcommand{\Re}{\operatorname{Re}}

\newcommand{\mo}{\sqrt{-1}}
\newcommand{\R}{\mathbb R}
\newcommand{\C}{\mathbb C}

\newcommand{\s}[1]{{\mathbb S}^{#1}}
\newcommand{\mcc}{\mathcal C}
\newcommand{\mci}{\mathcal I}
\newcommand{\mcj}{\mathcal J}
\newcommand{\trp}{\, {}^t \negthinspace}
\newcommand{\zb}{\bar{z}}
\newcommand{\dd}[2]{\frac{\partial {#1}}{\partial {#2}}}
\newcommand{\at}{\tilde{A}}

\newtheorem{thm}{Theorem}[section]

\newtheorem{cor}[thm]{Corollary}
\numberwithin{equation}{section}

\begin{document}

\allowdisplaybreaks

\renewcommand{\thefootnote}{$\star$}

\renewcommand{\PaperNumber}{068}

\FirstPageHeading

\ShortArticleName{Boundaries of Graphs of Harmonic Functions}

\ArticleName{Boundaries of Graphs of Harmonic Functions\footnote{This paper is a
contribution to the Special Issue ``\'Elie Cartan and Dif\/ferential Geometry''. The
full collection is available at
\href{http://www.emis.de/journals/SIGMA/Cartan.html}{http://www.emis.de/journals/SIGMA/Cartan.html}}}

\Author{Daniel FOX}

\AuthorNameForHeading{D.~Fox}

\Address{Mathematics Institute, University of Oxford, 24-29 St Giles', Oxford, OX1 3LB, UK}
\Email{\href{mailto:foxdanie@gmail.com}{foxdanie@gmail.com}}

\ArticleDates{Received March 06, 2009, in f\/inal form June 16, 2009;  Published online July 06, 2009}

\Abstract{Harmonic functions $u:\R^n \to \R^m$ are equivalent to integral manifolds of an exterior dif\/ferential system with independence condition $(M,\mci,\omega)$.  To this system one associates the space of conservation laws $\mcc$.  They provide necessary conditions for $g:\s{n-1} \to M$ to be the boundary of an integral submanifold.  We show that in a local sense these conditions are also suf\/f\/icient to guarantee the existence of an integral manifold with boundary $g(\s{n-1})$.  The proof uses standard linear elliptic theory to produce an integral manifold $G:D^n \to M$ and the completeness of the space of conservation laws to show that this candidate has $g(\s{n-1})$ as its boundary.  As a corollary we obtain a new elementary proof of the characterization of boundaries of holomorphic disks in $\C^m$ in the local case.}

\Classification{35J05; 35J25; 53B25}

\renewcommand{\thefootnote}{\arabic{footnote}}
\setcounter{footnote}{0}

\section{Introduction}
On $\C^m$ with complex coordinates $z^1,\ldots,z^m$, let
\[
\Omega^{(p,q)}=\left\lbrace f \ed z^{a_1} \w \cdots \w  \ed z^{a_p} \w \ed z^{b_1}  \w \cdots \w \ed z^{b_q}  \vert \;\; a_i,b_j \in \{ 1,\ldots,m\}, f \in C^{\infty}(\C^n,\C) \right\rbrace. \]
A holomorphic curve in $\C^m$ is a holomorphic map $\phi:X \to \C^m$, where $X$ is a Riemann surface.  This is equivalent to a \emph{real} $2$-dimensional surface $G:X \to \C^m$ for which
\begin{gather}\label{eq:HolomorphicCurve}
G^*\big( \Omega^{(2,0)}\big)=G^*\big( \Omega^{(0,2)}\big)=0.
\end{gather}
This equivalence can be demonstrated using the following argument.  Suppose that locally $G(X)$ can be written as a smooth graph $z^a=G^a(z^1,\zb^1)$.  Then $G^*(\ed z^a \w \ed z^1)=\dd{z^a}{\zb^1} \ed \zb^1 \w \ed z^1$ and thus~\eqref{eq:HolomorphicCurve} implies $\dd{z^a}{\zb^1} =0$.

The $1$-forms $\varphi$ satisfying $\ed \varphi \in \Omega^{(2,0)}\oplus \Omega^{(2,0)}$ provide moment conditions for the boundaries of holomorphic curves.  That is, for any holomorphic curve $G:X \to \C^m$ with boundary $g : \partial X \to \C^m$ we f\/ind
\begin{gather}\label{eq:Stokes}
\int_{\partial X} g^* \varphi = \int_{X}  G^*(\ed \varphi)=0.
\end{gather}
In particular, if $\varphi \in \Omega^{(1,0)}$ is holomorphic then $\ed \varphi   \in \Omega^{(2,0)}$ and so the integral of any holomorphic $(1,0)$-form around the boundary of a holomorphic curve is zero.  For more details see Example~4 of Section~1.1 of~\cite{Bryant1995}.

Using complex function theory Wermer \cite{Wermer1958} showed that, in the analytic category, the moment conditions provided by holomorphic $(1,0)$-forms are suf\/f\/icient to characterize the boundaries of holomorphic disks.  That is, if $D \subset \C$ is a domain with real analytic boundary $\partial D$, $g:\partial D \to \C^m$ is real analytic, and
\[
\int_{\partial D} g^* \varphi =0
\]
for all holomorphic $(1,0)$-forms on $\C^m$, then $g(\partial D)$ is the boundary of a holomorphic disk $G:D \to \C^m$.  This was generalized to higher dimensional domains for graphs by Bochner \cite{Bochner1943} and then a complete treatment was given by Harvey and Lawson \cite{Harvey1975,Harvey1977} using geometric measure theory.

We generalize this in another direction.  Using conservation laws (the analogue of the holomorphic $(1,0)$-forms) we characterize the boundaries of the graphs of $1$-jets of harmonic functions $\Delta u=0$ where $u:D \to \R^m$  has domain $D \subset \R^n$ with $C^2$ boundary (Theorem \ref{thm:bvforharmonic}).  From Theorem \ref{thm:bvforharmonic} we extract Corollary \ref{cor:bvforhol} which gives a new proof of the characterization of boundaries of embedded holomorphic disks.

We work using exterior dif\/ferential systems (EDS) and their characteristic cohomology, the relevant aspects of which we now review\footnote{For an introduction to EDS the reader might enjoy \cite{Ivey2003} or \cite{Bryant1991}.  The fundamental paper on characteristic cohomology is \cite{Bryant1995}.}.  To every partial dif\/ferential equation one can associate an EDS with independence condition. (In Section~\ref{sec:EDS} we do this for the Laplace equation.)  An EDS is a pair $(M,\mci)$ where $M$ is a manifold and $\mci \subset \Omega^*(M,\R)$ is a homogeneous dif\/ferential ideal.  An {\bf independence condition} is a totally decomposable nowhere vanishing form $\omega$ whose degree is the same as the dimension of the domain of the PDE.  The {\bf integral submanifolds}, that is, those submanifolds for which the ideal pulls back to be zero but $\omega$ pulls back to be nonzero, are equivalent to solutions of the PDE.  In this sense one can associate to every PDE a~submanifold geometry def\/ined by a non-degeneracy condition and the vanishing of dif\/ferential forms.

Each EDS $(M,\mci)$ (not necessarily with an independence condition) def\/ines cohomology groups $H^p(\Omega^*/\mci,\overline{\ed})$ on $M$ known as the characteristic cohomology \cite{Bryant1995}.  A certain graded piece of the characteristic cohomology constitutes the space of conservation laws.   Suppose that one is interes\-ted in $n$ dimensional integral submanifolds of $(M,\mci)$.  To each EDS one associates a positive integer ${{ l}}$, known as its characteristic number, that measures how overdetermined the EDS is (see Section~4.2 of~\cite{Bryant1995} for the def\/inition).   In \cite{Bryant1995} it is shown that over contractible open sets $H^{p}(\Omega^*/\mci,\overline{\ed})=0$ for $p<n-l$ when $\mci$ is involutive\footnote{See \cite{Ivey2003} or \cite{Bryant1991} for the def\/inition of involutivity.}.  The space of conservation laws $\mcc$ is def\/ined to be the f\/irst nontrivial cohomology group:
\begin{gather*}
\mcc=H^{n-l}(\Omega^*/\mci,\overline{\ed}).
\end{gather*}
For systems $(M,\mci)$ that arise from Lagrangians $l$ is equal to $1$.  In this case, a class in $\mcc$ is represented by a form $\varphi \in \Omega^{n-1}(M,\R)$ that is not in the ideal but for which $\ed \varphi \in \mci$.  From Stokes' theorem such forms lead to moment conditions on boundaries just as they did for holomorphic curves in equation~\eqref{eq:Stokes}.

\looseness=1
In Section~\ref{sec:EDS} we introduce the well known EDS with independence condition associated to the Laplace equation and a useful subspace of conservation laws.  In Section~\ref{sec:Boundaries} we use this set up to prove Theorem~\ref{thm:bvforharmonic}, which characterizes the boundaries of graphs of harmonic functions using the moment conditions arising from conservation laws.   In Section~\ref{sec:holomorphicdisks} we deduce Corollary~\ref{cor:bvforhol} which characterizes the boundaries of holomorphic disks that are graphs.

\section{The EDS for graphs of harmonic functions\\ and its conservation laws}\label{sec:EDS}

In this section we introduce the exterior dif\/ferential system and the space of conservation laws we will need. For the Laplace equation
\begin{gather}\label{eq:Laplace}
\Delta u =0
\end{gather}
where $u:\R^n \to \R^m  $ and $\Delta= (\dd{\;}{x^1})^2+\cdots +(\dd{\;}{x^n})^2$, we def\/ine
\[
M=J^1(\R^n,\R)=\R^n \times \R^m \times \R^m \otimes (\R^{n})^*
\]
to be the f\/irst jet space of maps from $\R^n$ to $\R^m$.  It has natural coordinates $(x^i,u^a,p^a_i)$.  The relevant dif\/ferential ideal is
\[
\mci=\langle \theta,\ed \theta, \psi \rangle,
\]
which is algebraically generated by the components of the vector valued dif\/ferential forms
\begin{gather*}
\theta=\ed u - p \ed x \in \Omega^1(M,\R^m),\\
\ed \theta = -\ed p \w \ed x \in \Omega^2(M,\R^m), \\
\psi=\ed p_i \w \ed x^{(i)} \in \Omega^n(M,\R^m),
\end{gather*}
where
\[
\ed x^{(i)}=*\ed x^i=(-1)^{i-1}\ed x^1 \w \cdots \w \widehat{\ed x^i} \w \cdots \w \ed x^n.
\]
Above $*$ is the Hodge star on $\R^n$ with respect to the standard f\/lat metric and volume form $\omega=\ed x^1 \w \cdots \w \ed x^n$.   We will use a mixture of index and matrix notation.  For example, $p \ed x$ is the $\R^m$-valued $1$-form with components $\sum_{i=1}^n p^a_i \ed x^i$ and in $\ed p_i \w \ed x^{(i)}$ the sum over $i$ is implicit.

The exterior dif\/ferential system $(M,\mci,\omega)$ for harmonic functions is involutive with characteristic number $l=1$.  Solutions to the Laplace equation \eqref{eq:Laplace} are equivalent to $n$-dimensional embedded submanifolds $G:X \to M$ such that $G^*(\mci)=0$ and $G^*\omega \neq 0$.   A $k$-dimensional submanifold $F:U \to M$ is def\/ined to be {\bf isotropic} if $F^*(\theta)=0$.

The general theory of characteristic cohomology of an exterior dif\/ferential system indicates that for the EDS associated to Laplace's equation the space of conservation laws is $\mcc=H^{n-1}(\Omega^*/\mci,\overline{\ed})$.  The short exact sequence
\[
0 \to \mci \to \Omega \to \Omega/\mci \to 0
\]
induces a long exact sequence in cohomology which, due to the vanishing $H^s_{dR}(M,\R)=0$ for $s>0$, produces the isomorphism
\[
\iota:\ H^{n-1}(\Omega^*/\mci,\overline{\ed}) \to H^n(\mci,\ed).
\]
The map is given by exterior dif\/ferentiation:  a class in $H^{n-1}(\Omega^*/\mci,\overline{\ed})$ is represented by a~dif\/fe\-ren\-tial form $\varphi \in \Omega^{n-1}(M)$ such that $\ed \varphi \in \mci$.  So if $[\varphi] \in H^{n-1}(\Omega^*/\mci,\overline{\ed})$ then $[\ed \varphi] \in H^n(\mci,\ed)$.  The class $[\ed \varphi]$ or its representative $\ed \varphi$ is referred to as the {\bf dif\/ferentiated} conservation law and~$[\varphi]$ or its representative $\varphi$ as the {\bf undif\/ferentiated} conservation law.  We will need the explicit form of dif\/ferentiated conservation laws and so turn to them now.

An element of $H^n(\mci,\ed)$ is represented by a closed $n$-form in $\mci$.  Any element in $\mci \cap \Omega^n(M,\R)$ is determined by an $\R^m$-valued $(n-1)$-form $\rho$, an $\R^m$-valued $(n-2)$-form $\sigma$, and an $\R^m$-valued function $H$, by the formula
\[
\Phi=\trp{\rho} \w \theta+\trp{\sigma} \w \ed \theta -\trp{H} \psi.
\]
We seek $\rho$, $\sigma$, $H$ that make $\Phi$ closed, but f\/irst we make a standard simplif\/ication.  Using the relation
\[
\Phi=\left(  \trp{\rho}-(-1)^{n-2}\ed \trp{\sigma} \right)  \w
 \theta -\trp{H} \psi  +(-1)^{n-2}\ed \left(\trp{\sigma} \w  \theta\right),
\]
and the fact that we are really only interested in the class  $[\Phi] \in H^n(\mci,\ed)$, we see that, for any class in $H^n(\mci,\ed)$, we can always f\/ind a representative for which $\sigma=0$.

The following special set of conservation laws will be suf\/f\/icient for studying the boundaries of integral manifolds that satisfy the independence condition.  If $H:\R^n \to \R^m$ is a harmonic function and $\rho=(-1)^{n}*\ed H $, then
\begin{gather}\label{eq:NormalForm}
\Phi=\trp{\rho} \w \theta +\trp{H} \psi
\end{gather}
is closed and represents a class in $H^n(\mci,\ed)$.

Conservation laws are a natural source of moment conditions.  Let $\varphi$ be an undif\/ferentiated conservation law.  By stokes theorem
\begin{gather*}
\int_{\partial D} g^* \varphi = \int_{D}  G^*(\ed \varphi)
\end{gather*}
for any $G:D^n \to M$ with $g=G_{\vert_{\partial D}}$.  Thus if $G(D)$ is integral
\begin{gather*}
\int_{\partial D} g^* \varphi =0.
\end{gather*}
In the next section we show that in a local sense the moment conditions coming from conservation laws are complete for the harmonic function system.

\section{Boundaries of graphs of harmonic functions}\label{sec:Boundaries}
Let $\pi:M \to \R^n\times \R^m$ be the standard projection $(x,u,p) \mapsto (x,u)$.
\begin{thm}\label{thm:bvforharmonic}
Let $g:\s{n-1} \to M$ be a $C^2$ isotropic submanifold such that $x \circ g:\s{n-1} \to \R^n$ is an embedding and $x \circ g (\s{n-1})$ is the boundary of a domain $D \subset \R^n$.  Then there exists $G:D \to M$ such that
\begin{gather*}
G^*(\mci)=0,\qquad
G(\partial D)=g\big(\s{n-1}\big)
\end{gather*}
if and only if
\begin{gather}\label{eq:MomentConditions}
\int_{g(\s{n-1})}\varphi=0 \qquad \forall \, \varphi \in \mcc.
\end{gather}
\end{thm}

The proof relies on the standard theory of linear elliptic PDE to produce an integral submanifold and then uses the moment conditions arising from conservation laws to show that it has the desired boundary.  We use the following existence and uniqueness result \cite{Gilbarg2001}.
\begin{thm}\label{thm:dirichlet}
Let $D \subset \R^n$ be a bounded domain with $C^2$ boundary $\partial D$ and let $v:\partial D \to \R^m$ be continuous.  Then there is a unique function $V:D \to \R^m$ satisfying $V_{{\vert}_{\partial D}}=v$ and $\Delta V=0$.
\end{thm}

\begin{proof}[Proof of Theorem~\textup{\ref{thm:bvforharmonic}}]
Let $v=u_{|_{\partial D}}$.  Then by Theorem \ref{thm:dirichlet} there exists a unique smooth function $V:D \to \R^m$ such that $\Delta V=0$ and $V_{| \partial D}=v$.  Let $J^1(V):D \to M$ be the $1$-jet of $V$, $J^1(V)(x)=(x,V(x), \nabla V(x))$.  By construction $\pi \circ g(\s{1})=\pi \circ J^1(V)(\partial D)$.  We will now show that{\samepage
\[
J^1(V)(\partial D)=g\big(\s{n-1}\big).
\]
Once this is accomplished, $G=J^1(V)$ is the desired solution.}

Let $\tilde{g}=J^1(V)_{|\partial D} \circ x \circ g:\s{n-1} \to M$.  Then $\tilde{g}, g:\s{n-1} \to M$ are isotropic submanifolds that both satisfy the moment conditions \eqref{eq:MomentConditions} and $\pi \circ g= \pi \circ \tilde{g}$.  Write
\begin{gather*}
g(s)=\big(x(s),u(s),A(s)\big),\qquad
\tilde{g}(s)=\big(x(s),u(s),\at(s)\big),
\end{gather*}
where $\at=\nabla V_{|\partial D} \circ x \circ g$ and $s \in \s{n-1}$.   For $i=1 ,\ldots, n-1$ let $s^i$ be local coordinates on $\s{n-1}$. The fact that $g:\s{n-1} \to M$ is isotropic implies that
\[
0=g^*\theta^a=\left(\dd{u^a}{s^i}-A^a_j \dd{x^j}{s^i}\right)\ed s^i
\]
so that
\[
 \dd{u^a}{s^i}=A^a_j \dd{x^j}{s^i}.
\]
Similarly we f\/ind that
\[
\dd{u^a}{s^i}=\at^a_j \dd{x^j}{s^i}.
\]
Let $N$ be the outward unit normal of $x \circ g(\s{n-1}) \subset \R^n$ so that
\[
0=N_i \dd{x^i}{s^j}.
\]
We can now decompose $\at(s)=A(s)+\zeta(s) \trp{N}(s) \in \R^m \otimes (\R^n)^*$ for some function $\zeta:\s{n-1} \to \R^m$.

Let
\[
\chi: \ [0,1]  \times \s{n-1} \to M
\]
be the smooth map given in coordinates by
\[
\chi(r,s)=(x(s),u(s),p(r,s)),
\]
where
\[
p(r,s)=A(s)+(1-r)\zeta(s)\trp{N(s)}.
\]
The cylinder $\chi(\s{n-1} \times I)$ has the following properties:
\begin{itemize}\itemsep=0pt
\item{$\chi(0,s)=\tilde{g}(s)$};
\item{$\chi(1,s)=g(s)$};
\item{$\chi^*\theta^a=0$};
\item{$\chi(\s{n-1} \times I) \subset \pi^{-1}\big( \pi  \circ g (\s{n-1})\big)$}.
\end{itemize}
Therefore $\partial( \chi(\s{n-1} \times I) )=g(\s{n-1})-\tilde{g}(\s{n-1})$ and because both $g$ and $\tilde{g}$ satisfy the moment conditions induced by conservation laws $\Phi=\ed \varphi$,
\begin{gather}\label{eq:PhiVanishes}
\int_{\chi(\s{n-1} \times I)}\Phi=\int_{g(\s{n-1})}\varphi-\int_{\tilde{g}(\s{n-1})}\varphi=0.
\end{gather}
Now assume that $\Phi$ is of the type specif\/ied in \eqref{eq:NormalForm}.  Because $\chi$ is contact,
\[
\chi^*(\Phi)=\chi^*(\trp{H} \; \psi)
\]
and we calculate that
\begin{gather*}
\chi^*(\psi^a)= \chi^*(\ed p^a_i  \w \ed x^{(i)})   =  -\zeta^a N_i \ed r  \w \chi^*(\ed x^{(i)})
 =- \zeta^a \ed r \w \chi^*(N \lhk \omega)=- \zeta^a \ed r \w g^*(N \lhk \omega).\!
\end{gather*}
The def\/inition of $\chi$ also implies that $\chi^*(H)=g^*(H)$.  By Theorem \ref{thm:dirichlet}, the harmonic function $H:B^n \to \R^m$, where $B^n$ is the closed ball with boundary $\s{n-1}$, is uniquely determined by choosing an arbitrary continuous function, $h:\s{n-1} \to \R^m$, and specifying that  $g^*(H)=h$.

We now calculate
\begin{gather*}
\int_{\chi(\s{n-1} \times I)}\Phi = \int_{I \times \s{n-1}} \chi^*(\trp{H}\psi)
=-\int_{I \times \s{n-1}} \trp{h} \zeta \ed r \w g^*(N \lhk \omega)\\
\phantom{\int_{\chi(\s{n-1} \times I)}\Phi}{} =-\int_{\s{n-1}}\trp{h}\zeta \cdot \left( \int_{0}^1 \ed r \right)  g^*(N \lhk \omega) =-\int_{\s{n-1}} \trp{h}\zeta  g^*(N \lhk \omega).
\end{gather*}
Using \eqref{eq:PhiVanishes} this implies that
\[
0=\int_{\s{n-1}} \trp{h} \zeta g^*(N \lhk \omega)
\]
for all continuous functions $h:\s{n-1} \to \R^m$.  The $(n-1)$-form $g^*(N \lhk \omega)$ is the induced volume form on $x \circ g:\s{n-1} \to \R^n$.  Therefore, using the assumption that $x \circ g$ is an embedding, we can conclude that $\zeta=0$. This implies that $A=\at$ and thus $\tilde{g}=g$.
\end{proof}

\section{Boundaries of holomorphic disks}\label{sec:holomorphicdisks}

On $\C^{m}$ let $\hat \mcj=\langle \Omega^{2,0} \oplus \Omega^{0,2} \rangle$.  Then a holomorphic curve is a real surface that is an integral manifold for $\hat \mcj$.  Using distinct approaches, Wermer \cite{Wermer1958} and Harvey and Lawson \cite{Harvey1975} prove
\begin{thm}
Let $Y \subset \C^m$ be a compact, connected, oriented submanifold of dimension one and of class $C^2$.  Suppose that
\[
\int_Y \varphi=0
\]
for all holomorphic $1$-forms $\varphi \in \Omega^{(1,0)}$.  Then there exists an irreducible holomorphic curve $X \in \C^m \setminus Y$ such that $\partial X=Y$.
\end{thm}

Wermer provided the f\/irst such result using complex function theory. Harvey and Lawson actually prove a much stronger result that characterizes the boundaries of complex submanifolds in which the boundary may have multiple connected components.    We can deduce a local version of this from Theorem \ref{thm:bvforharmonic}.

\begin{cor}\label{cor:bvforhol}
Let $g:\s{1} \to \C^{m}$ be a $C^2$ embedded curve for which there exists a projection to a complex line $\zeta:\C^{m} \to \C$ such that $ \zeta \circ g(\s{1})$ is the boundary of a domain $D \subset \C$.  Then there exists a holomorphic map $G:D \to \C^{m}$ such that
\[
G({\partial D})=g\big(\s{1}\big)
\]
if and only if
\begin{gather}
\int_{g(\s{1})}\varphi=0
\end{gather}
for all holomorphic $1$-forms $\varphi$ on $\C^m$.
\end{cor}

\begin{proof}
We make an integrable extension and then rely on Theorem \ref{thm:bvforharmonic}.  Let $(M,\mci,\omega)$ be the system for harmonic functions $u:\R^2 \to \R^{m-1}$, so that
\[
M=\R^2 \times \R^{m-1} \times  \R^{m-1} \otimes \big(\R^2\big)^*.
\]
Let
\[
\zeta:M \to \R^2 \times \R^{m-1} \otimes \big(\R^2\big)^*
\]
be the standard projection and identify the image with $\C^m$ by def\/ining the holomorphic coordinates $z^a=p^a_1+\mo p^a_2$ for $a=1, \ldots, m-1$ and  $z^m=x^1-\mo x^2$. We def\/ine the dif\/ferential ideal $\mcj=\langle \ed z^a \w \ed z^m \rangle$ on $\C^m$ and let $\Omega=-\frac{\mo}{2}\ed z^m \w \ed \zb^m$ def\/ine an independence condition. The integral manifolds of $(\C^m,\mcj,\Omega)$ are holomorphic disks that can be graphed as functions of $z^m$.   It is readily checked that $\zeta^*(\ed z^a \w \ed z^m) = \ed \theta^a - \mo \psi^a$ and that $\zeta^*(\Omega)=\omega$.  Therefore the projection $\zeta:M \to \C^{m}$ makes $(M,\,\mci,\omega)$ an integrable extension of $(\C^m,\mcj,\Omega)$:  that is, $\mci$ is {\it algebraically} generated by $\zeta^*\mcj$ and the $1$-forms $\theta^a$, and the independence conditions are compatible.

First we must show that when the conservation laws of $(\C^m,\mcj,\Omega)$ are pulled back using $\zeta$, they surject onto the special class of conservation laws used in the proof of Theorem~\ref{thm:bvforharmonic}.  Then we must show that we can lift the supposed boundary $g:\s{1} \to \C^m$ to $f:\s{1}\to M$ so that $f(\s{1})$ satisf\/ies all of the moment conditions for $(M,\mci,\omega)$.  For the f\/irst part we must show that for any harmonic function $H:D \to \R^{m-1}$ and $\Phi$ def\/ined from $H$ as in \eqref{eq:NormalForm}, there is a conservation law for $\mcj$ that pulls back under $\zeta$ to give the same class as $[\Phi] \in H^2(\mci,\ed)$.  To see this, let $H$ be the desired harmonic function and let $K:D \to \R^{m-1}$ be its harmonic conjugate, so that $K^a + \mo H^a$ is a holomorphic function of $z^m$. Let
\[
\Upsilon=\big(K^a + \mo H^a\big) \ed z^a \w \ed z^m \in \mcj.
 \]
Then $\ed \Upsilon =0$ and $\Upsilon$ is a dif\/ferentiated conservation law for $\mcj$.  When pulled up to $M$ we f\/ind
\[
\zeta^*(\Re(\Upsilon))=K^a \ed \theta^a + H^a \psi^a.
\]
We can rewrite this as
\[
\zeta^*(\Re(\Upsilon))=\trp \rho  \w \theta + \trp H \psi + \ed \big(\trp H \theta\big),
\]
where $\rho=* \ed H$, which is consistent with \eqref{eq:NormalForm}.  Since $\ed \Upsilon=0$, the form $\trp \rho \w \theta+\trp H \psi$ is also closed and therefore a dif\/ferentiated conservation law of the form \eqref{eq:NormalForm} with the desired harmonic function~$H$.  Therefore if $g:\s{1} \to \C^m$ satisf\/ies all of the moment conditions from the conservation laws of~$\mcj$, then an appropriate lift to $M$ will satisfy all of the moment conditions for $\mci$ that are needed to apply Theorem~\ref{thm:bvforharmonic}.

Now we turn to the lift.  Any lift $f:\s{1} \to M$ of $g$ is def\/ined by choosing a map $u:\s{1} \to \R^{m-1}$.  For the lift to be contact we must have
\[
0=f^*(\theta^a)=\dd{u^a}{s}-\left(p^a_1\dd{x^1}{s}+p^a_2\dd{x^2}{s}\right),
\]
so def\/ine
\[
u^a(s)= \int_0^s \left( p^a_1(t)\dd{x^1}{t}+p^a_2(t)\dd{x^2}{t} \right) \ed t.
\]
This is a periodic function since $g$ satisf\/ies the moment condition
\[
\int_{g(\s{1})}\big(p^a_1\ed x^1+p^a_2 \ed x^2\big)=0.
\]
Now by Theorem \ref{thm:bvforharmonic} there exists $F:D^2 \to M$ such that $F_{| \partial D}=f$ and $F^*(\mci)=0$.  Then $G=\zeta \circ F:D^2 \to \C^m$ is the desired holomorphic disk.
\end{proof}

\subsection*{Acknowledgements}
I would like to thank Dominic Joyce and Yinan Song for useful conversations.  This work was carried out with the support of the National Science Foundation grant OISE-0502241.

\pdfbookmark[1]{References}{ref}
\LastPageEnding

\end{document}